\documentclass[12 pt,english,letter, french]{article}
\usepackage{amsmath}
\usepackage{amssymb}
\usepackage{graphicx}
\usepackage[ansinew]{inputenc}

\def\R{\mathbb R}

\def\epsilon{\varepsilon}

\def\epsilon{\varepsilon}
\def\ds{\displaystyle}

\numberwithin{equation}{section}

\newtheorem{thm}{Theorem}[section]

\newtheorem{notation}[thm]{Notation}

\title{The non-monotonicity of  the KPP speed with respect to diffusion in the presence of a shear flow}
\author{Mohammad El Smaily \footnote{Email: elsmaily@andrew.cmu.edu} \thanks{The author is indebted to the CMU-Portugal program and ``Center for Nonlinear Analysis'' for their support.}\\
\footnotesize{ Department of Mathematical Sciences and Center for Nonlinear Analysis,}\\
\footnotesize{Carnegie Mellon University, Wean Hall, }\footnotesize{5000 Forbes Avenue, Pittsburgh, PA, 15213, USA}\\
}

\date{February 28, 2011}
\begin{document}
\maketitle

\begin{abstract} In this paper, we prove via counterexamples that adding an advection term of the form Shear flow (whose streamlines are parallel to the direction of propagation) to a reaction-diffusion equation will be an enough heterogeneity to spoil the increasing behavior of the KPP speed of propagation with respect to diffusion. The non-monotonicity of the speed with respect to diffusion will occur even when the reaction term and the diffusion matrices are considered homogeneous (do not depend on space variables). For the sake of completeness, we announce our results in a setting which allows domains with periodic perforations that may or may not be equal to the whole space $\mathbb{R}^N.$

\vskip 0.5cm

\noindent{\it Keywords:} \rm Traveling fronts, reaction-diffusion, monotonicity with respect to diffusion, KPP minimal speed,  shear flows, principal eigenvalue.\hfill\break

\end{abstract}

\section{Introduction and motivation}\label{intro}

Pulsating traveling fronts are particular solutions of heterogeneous Reaction-Advection-Diffusion equations that describe propagation phenomena for combustion models, evolution of epidemics, population dynamics 
and many other phenomena. This paper is dedicated to answer the question whether the minimal speed of propagation (KPP nonlinearity) is monotone with respect to the diffusion coefficient 
or not. We answer this question in the case where an advection term is present. This work is also a continuation of the paper \cite{elsmaily1} which dealt with the dependence and asymptotic behaviors of the minimal speed of propagation on the coefficients of the reaction-diffusion equation in heterogenous media. 

To explain the question mathematically we start by recalling the definition of the minimal speed of propagation and the changes that happen while passing from a homogeneous to a heterogeneous setting.   
\subsubsection{Homogeneous framework}\label{subs1}
 In the article \cite{KPP} by Kolmogorov, Petrovsky and Piskunov, the setting was ``homogeneous'' in the following sense. The equation considered in \cite{KPP} was 
\begin{equation}\label{hom_reaction_diffusion}
u_t=\Delta u +f(u)\quad \hbox{for all }\quad(t,x)\in\mathbb{R}\times\mathbb{R}^N,
\end{equation}
where $0\leq u=u(t,x)\leq1$ was a function defined over $\R\times\R^N.$
The reaction term $f$ was assumed to be of KPP type. That is, $f:[0,1]\mapsto\R$ such that
\begin{equation}\label{kpp_hom_1}
\begin{array}{c}
 f(u)>0 \hbox{ for all } u\in (0,1),\quad  f'(0)>0,
\end{array}
\end{equation}
 together with the sub-linearity condition (which is usually called the ``KPP condition'') over the interval $[0,1]$
 \begin{equation}\label{kpp_hom_2}
\forall u\in[0,1],\;f(u)\leq f'(0)u.
\end{equation}
One can see that, in the equation (\ref{hom_reaction_diffusion}), the diffusion and reaction terms do not depend on space and time variables and there is no advection term of the form $q(x)\cdot \nabla u$. Fixing a 
unitary direction $e\in\mathbb{\R}^N,$ traveling fronts in the direction of $-e$ and with a speed $c,$ in such a homogeneous setting, are solutions $u(t,x)=\phi(x\cdot e+ct)=\phi(s)$ which satisfy the limiting 
conditions $\phi(-\infty)=0$ and $\phi(+\infty)=1.$ Kolmogorov \textit{et al.} \cite{KPP} proved 
that when the reaction term is of KPP type (i.e.~satisfying the conditions (\ref{kpp_hom_1}) and (\ref{kpp_hom_2})), there exists a minimal speed  $c^*=2\sqrt{f'(0)}$ such that a traveling front propagating in the 
direction of $-e$ with a speed $c$ exists if and only if $c\geq c^*=2\sqrt{f'(0)}.$
If we look at the equation 
\begin{equation}\label{homeq1}
u_t=D\Delta u+f(u),
\end{equation}
for some positive constant $D,$ then a simple rescaling yields that the minimal KPP speed of (\ref{homeq1}) is given by $c^*_{D}=2\sqrt{D}\sqrt{f'(0)}.$ One can then notice that, in a homogeneous setting, the minimal speed is increasing with respect to the diffusion factor (the map $D\mapsto c^*_{D}$ is increasing over $(0,+\infty)$). 

 \subsubsection{Heterogeneous framework, notations and setting}
  The result of Kolmogorov et al. \cite{KPP} has been generalized to media with periodic spatially dependent coefficients (see \cite{bh1}, \cite{shigesada Kawsaki 2}, \cite{weinberger}, \cite{JXin} for example) and to 
settings with space-time dependent coefficients (see \cite{bhnadin, nadin1, nadin2}). We recall here some results which are very relevant to this present work and  we also introduce the setting that will be considered.  
 The equation that we will consider is of the type
\begin{equation}\label{heteq}
u_t(t,z)=\nabla_z\cdot(A(z)\nabla_z u)+q(z)\cdot \nabla_z u +f(u),\quad t\in\R, z\in \Omega,
\end{equation} 
where the domain $\Omega:=\R\times \omega$ is an unbounded $C^3$ open connected subset of $\R^N=\R\times \R^{N-1}$ ($N\geq 1$). We represent each $z\in\Omega$ as 
$$z=(x,y),~~x\in \R, ~ y\in \omega.$$
 In order to cover all possible cases, we will now give a generic description of $\omega$ which appears in $\Omega$ when the space dimension is $N>1.$ The set $\omega$ is assumed to have the following 
structure. There exists an integer $d\in\{0,\ldots,N-1\}$ and $L_1>0,\ldots, L_d>0$ so that an element $y\in\omega$ is represented
 by $(y_1,y_2)\in \R^{d}\times\R^{N-1-d}$ and  

\begin{eqnarray}\label{comega}
    \left\{
      \begin{array}{l}
        \exists\,R\geq0\,;\forall\,(y_1,y_2)\in\omega,\,|y_2|\,\leq\,R, \vspace{3pt} \\
        \forall\,(k_1,\cdots,k_d)\;\in\;L_1\mathbb{Z}\times\cdots\,\times L_d\mathbb{Z},
        \quad\displaystyle{\omega=\omega+\sum^{d}_{k=1}k_ie_i},
      \end{array}
    \right.
\end{eqnarray}
where $(e_i)_{1\leq i\leq N}$ is the canonical basis of $\mathbb{R}^N.$

\begin{itemize}
\item Notice that the case $d=0$ corresponds to ``$\omega$ is bounded'' and hence $\Omega=\R\times\omega$ is an infinite cylinder of section $\omega.$
\item In the cases where $d\geq1,$ the set $\omega$ is unbounded. 
\end{itemize}
As particular examples of $\Omega,$ one can have the whole space $\R^N=\R\times\R^{N-1}$ in which case $d=N-1.$ One can also have the whole space $\R^N$ except a periodic array of holes (periodic perforations). In such case, $d$ is also equal to $N-1.$  For $N=2,$  one has $d\in\{0,1\}.$ The case where $d=1$ (with $N=2$) means that $\omega$ is unbounded connected and satisfies (\ref{comega}) and thus $\Omega=\R^2$. For $N=3,$ $d$ can be $0,1,$ or 2. The case where $d=1$ corresponds to $\omega\subset \R^2$ bounded in one direction and unbounded in the other one. In the unbounded direction, $\omega$  has to be 
periodic with a period that we denoted by $L_1$ in the general setting above. The case where $d=2$ corresponds to $d=N-1$, and thus, has been discussed above.    

 We mention that in all cases (even when $d=0$), the domain $\Omega:=\R\times\omega$ has a periodicity cell which we denote by $C:=[0,1]\times C_w$ where 
 \begin{equation}\label{cell}
\left\{\begin{array}{ll}
 C_\omega=&\omega, \text{ when } d=0, \vspace{4 pt}\\
C_\omega=&\{(y_1,y_2)=(y_1^1,\ldots,y_1^d,y_2) \in \omega,  y_1^1\in[0, L_1]\ldots,y_1^d\in [0, L_d]\} \text{ otherwise.}
 \end{array} \right.
  \end{equation}

  Concerning the reaction term,  we will only deal with nonlinearities $f$ depending on $u$ in order  achieve the proof of non-monotonicity of the speed with respect to diffusion. Precisely, $f$ is of the homogeneous 
KPP type (\ref{kpp_hom_1}-\ref{kpp_hom_2}).

  The advection term, in this paper, is assumed to be a shear flow. That is a vector field $q(x,y)=(q_1(y),0,\ldots,0)$ of class 
$C^{1,\alpha}(\overline{\Omega})$  ($\alpha>0$). The advection is assumed to satisfy  
\begin{equation}\label{cq}
       \left\{ \begin{array}{ll}
        q_1\hbox{ is $(L_1,\ldots,L_d)-$periodic with respect to }y \text{ (whenever  $d\geq1$)},\vspace{4 pt}\\
        \displaystyle{\int_{C}q_1(y)dxdy =0}\hbox{.}
\end{array}\right.
\end{equation}
Obviously the above vector field satisfies 
$$ \nabla_{x,y}\cdot q=0\hbox{ in }\overline{\Omega} ~~\hbox{ and }~~q\cdot\nu=0 \hbox{ on }\partial\Omega,$$ where $\nu$ is the outward normal on $\partial \Omega$ and it is given by $\nu=(0,\nu_{\omega})$ 
where $\nu_\omega$ is the outward normal on $\partial\omega.$ 

Let us now describe a wide class of diffusion coefficients for which the existence of traveling fronts and minimal speed of propagation holds according to \cite{bh1}. The diffusion matrix $A(x,y)=A(x,y_1,y_2)=(A_{ij}
(x,y))_{1\leq i,j\leq
 N}$ is a symmetric $C^{2,\alpha}(\,\overline{\Omega}\,)$ (with $\alpha>0$) matrix field satisfying
\begin{eqnarray}\label{cA}
    \left\{
      \begin{array}{l}
        A\; \hbox{is $(1,L_1,\cdots,L_d)-$periodic with respect to}\;(x,y_1),  \vspace{4 pt}\\
        \exists\,0<\alpha_1\leq\alpha_2;\forall(x,y)\;\in\;\Omega,\forall\,\xi\,\in\,\mathbb{R}^N, \vspace{4 pt}\\
        \displaystyle{\alpha_1|\xi|^2 \;\leq\;\sum_{1\leq i,j\leq N} \,A_{ij}(x,y)\xi_i\xi_j\,\;\leq\alpha_2|\xi|^2.}
      \end{array}
    \right.
\end{eqnarray}

In the above setting where $\Omega=\R\times\omega$ satisfies (\ref{comega}),  $q$ is of the type (\ref{cq}) and $A$ satisfies (\ref{cA}), one can rewrite (\ref{heteq}) as 
\begin{equation}\label{heteq1}
    u_t(t,x,y)=\nabla\cdot(A(x,y)\nabla u)+q_{_1}(y)\partial_{x}u+f(u) \hbox{ in }\R\times \R\times\omega,
\end{equation}
together with the boundary condition (in the cases where $\partial \omega\neq\emptyset$)
\begin{equation}\label{boundary}
 \nu(x,y)\cdot A\nabla_{x,y}u(t,x,y)=0,~(t,x,y)\in\mathbb{R}\times\mathbb{R}\times\partial\omega.
\end{equation}

Now in this non-homogenous setting,  we set $e=(1,0,\ldots,0)\in\R^{N}$ as the direction of propagation. We recall, from \cite{bh1, weinberger, JXin}, that a pulsating traveling front in the direction of $-e$ that 
propagates with a speed $c\in\R$ is a classical  solution $u=u(t,x,y):=\phi(x+ct,x,y)$  that connects $0$ to $1$ as follows  $$\lim_{x\rightarrow -\infty}u(t,x,y)=0\hbox{ and }\lim_{x\rightarrow+\infty}u(t,x,y)=1$$ 
(locally in $t$ and uniformly in $(x,y)$), and satisfies $\displaystyle{u(t+\frac{1}{c},x,y_1+k,y_2)=u(t,x+1,y_1,y_2)}$ for any $k\in L_1\mathbb{Z}\cdots\times L_d\mathbb{Z}.$ In terms of $\phi,$ the latter means that $\phi \hbox{ is 1-periodic in }x$ and $(L_1,\cdots,L_d)-$periodic in $y_1.$
 Under the assumptions (\ref{comega}-\ref{cq}-\ref{cA}) on $\Omega,$ $q$ and $A$, and having a KPP nonlinearity 
(\ref{kpp_hom_1}-\ref{kpp_hom_2}) one knows from \cite{bh1} that there exists a minimal speed $c^*:=c^*_{A,\Omega, q,f}(e)>0$ so that a pulsating traveling front with a speed $c$ exists if and only if $c\geq c^*_
{\Omega,A,q,f}(e).$
The minimal speed $c^*$ has been expressed in terms of the coefficients of the reaction-advection-diffusion problem via a variational formula which involves eigenvalue problems in \cite{BHN1} and \cite
{weinberger}. This formula has been used in many works to study various asymptotic and homogenization regimes of pulsating traveling fronts (see for eg.  \cite{elsmaily1}, \cite{ek1}, \cite{llm1}, \cite{nadin}, \cite{Zlatos Ryzhik}, \cite{zlatos}). 
 
  \section{The non-monotonicity of the minimal speed with respect to diffusion (in the presence of an advection field)}\label{counter_examples}

 Berestycki, Hamel, and Nadirashvili \cite{BHN1}
proved (in part $2$ of Theorem $1.10$) that, having any periodic
domain $\Omega\subseteq\mathbb{R}^{N}$,
$q\equiv0$ and a constant growth rate $f'_u(x,y,0)$ (which holds in the case $f=f(u)$), the map $\displaystyle{\beta\mapsto
c_{\Omega,\beta A,0,f}^*(e)}$ is increasing in $\beta>0.$ Nadin \cite{nadin} has proved, in Theorem 2.5, that the same result still holds even when $f'_u(x,y,0)$  depends on the spatial variables $(x,y)$. However,  we notice in both results of \cite{BHN1, nadin}  the absence of an advection term. 

 Furthermore,  in Theorem 6.1 of \cite{elsmaily1}, where
$\Omega=\mathbb{R}\times\omega,$ $q$ is a shear flow of the form
$q(x,y)=(q_1(y),0,\ldots,0)$ in $\overline{\Omega},$ while $A$  and
$f$ satisfy certain assumptions which are valid in our present setting, it was proved that the map
$\displaystyle{\beta\mapsto\frac{c_{\Omega,\beta
A,\sqrt{\beta}\,q,f}^*(e)}{\sqrt{\beta}}}$ is decreasing with
respect to $\beta>0$ in both cases $q_1\not\equiv0$ or $q_1\equiv0$
over $\omega.$ Having those results in \cite{BHN1}, \cite{elsmaily1} and \cite{nadin}, there arise naturally the
following two questions.
\begin{itemize}
   \item First: Do we still have the increasing behavior of the minimal
speed with respect to the diffusion factor $\beta$ in the
presence of an advection, even if the nonlinearity is homogenous?
\item Second:  Owing to  Theorem 1.1 in
\cite{BHN1}, the map $A\mapsto c_{\Omega, A,q,f}^*(e),$
 where $A$ varies in the ordered family \footnote{We say that $A=A(x,y)\leq B=B(x,y)$ if and only if for each $(x,y)\in\Omega$ and for each $\xi\in\mathbb{R}^{N},$ we have
$\xi \cdot A(x,y)\xi \leq \xi \cdot B(x,y)\xi.$ Also, we say that $A<B$ if and only if
for each $(x,y)\in\Omega$ and for each $\xi \in\mathbb{R}^{N}\setminus\{0\},$ we have
$\xi \cdot A(x,y)\xi<\xi \cdot B(x,y)\xi.$} of positive definite matrices satisfying
(\ref{cA}) is well defined (provided that $\Omega,$ $q$ and
$f$ satisfy (\ref{comega}), (\ref{cq}) and (\ref{kpp_hom_1}-\ref{kpp_hom_2})). We investigate the variation of the minimal speed of
propagation with respect to the matrix of diffusion. More
precisely, if $A$ and $B$ are two positive definite
matrices satisfying (\ref{cA}) and if $A<B,$ do we still have
$c_{\Omega, A,q,f}^*(e)\leq c_{\Omega, B,q,f}^*(e)$ whenever $q\not\equiv 0$?
\end{itemize}

As a matter of fact, the presence of a shear flow will change the monotone behavior of the minimal speed with respect to diffusion. We
prove in Subsection \ref{counter_to_first_question} that the answer is negative in general even when the diffusion matrices
$A$ and $B$ are proportional to the identity matrix and the nonlinearity is homogenous. We give a counterexample when the advection is large (up to a scaling, this gives a counterexample to the monotonicity of
the speed with respect to diffusion with a prefixed term and a small reaction). In Subsection \ref{counter to second question},  we prove that,  for a fixed nonzero shear flow and a fixed reaction, the answer to the second question is negative in general for matrices $A\leq B$ which are not equal up to a positive scalar. We mention that the second result cannot be seen as a trivial consequence of the first one by using the argument that ``a small perturbation of the strict inequality (\ref{delicate}) remains a strict inequality and one can then perturb the diffusion matrices'' (see the more precise details at the beginning of Subsection \ref{counter to second question}). Lastly, we mention that our counterexamples apply when the diffusion matrices do not depend on the variable $x$. One knows that, in such case, the traveling fronts of the reaction-advection-diffusion problem (\ref{heteq1}) will have the form $\phi(x+ct,y)$ and they will be invariant with respect to the frame moving in the direction $-x$ (or $-e$). Moreover, we know that in a homogenous setting the speed is monotone with respect to diffusion factors (see Section \ref{subs1}). Thus, the closer we are to a homogeneous setting the harder it is to prove non-monotonicity of the KPP speed with respect to diffusion.  Indeed, the setting where we construct our counterexamples can be  taken very close to the homogeneous one with only one heterogeneity coming from the advection. This gives sharpness to our results.  

\subsection{A counterexample devoted to answer the first question}\label{counter_to_first_question}
In this subsection, we will show a reaction-advection-diffusion problem whose diffusion matrix varies in the subfamily of positive definite matrices $\displaystyle{PD_{Id}}=\{\beta\,Id,~~\beta>0\}$ (where $Id$ stands for the $N\times N$ identity matrix) while  a shear flow 
will make the Part $2$ of Theorem $1.10$ \cite {BHN1} not valid anymore, even if the nonlinearity $f$ is taken homogeneous. Let's first announce the following result.

\begin{thm}\label{proportional}
Let $e=(1,0,\ldots,0)\in\mathbb{R}^{N},$ $\Omega=\mathbb{R}\times\omega\subseteq\mathbb{R}^{N}$ satisfying (\ref{comega})
where $\omega$ may or may not be bounded.  Assume that the
nonlinearity $f=f(u)$ is a homogenous ``KPP'' nonlinearity satisfying (\ref{kpp_hom_1}-\ref{kpp_hom_2}) and let's start with a shear flow $q(x,y)=(q_{1}(y),0,\ldots,0)$ defined over
 $\Omega$ such that $q_{1}\not\equiv0$ (satisfying (\ref{cq}) when $d \geq1$)  and so that
\begin{equation}\label{smaily}
0<2\sqrt{f'(0)}+\delta<\max_{y\in\overline{\omega}}\left(q_{1}(y)\right)-\delta,~~\hbox{for some $\delta>0.$}
\end{equation}
Then, there exists $M_1:=M_1(\delta)>0$ for which there corresponds $0<\epsilon_1:=\epsilon_1(\delta,M_1) < M_1$ such that
\begin{equation}\label{delicate}
\forall 0<\epsilon \leq \ds \epsilon_1< M_1,~~c_{\textstyle{\Omega,\epsilon\,\textstyle{Id},\sqrt{M_1}\,q, f}}^*(e)>c_{\textstyle{\Omega,M_1\textstyle{Id},\sqrt{M_1}q, f}}^*(e).
\end{equation} 
\end{thm}

\noindent \textbf{Proof.} Before going further into details, we recall that the variational formula of the KPP minimal  speed (see \cite{BHN1} for example) yields the continuity of $(\kappa,\rho,\mu)\mapsto c^*_{\Omega,\kappa A,\rho q,\mu f}(e)$
with respect $\kappa>0,$ $\mu>0$ and $\rho \in \R.$ This continuity with respect to the factors of the reaction, diffusion and advection will be useful to construct our proof. Due to the presence of several parameters in the problem, we proceed in simple steps.

\underline{Step 1.} In the author's work \cite{elsmaily1}, an asymptotic regime for the speed $c^*$ within large diffusions $M A$ together with advection fields of the form $M^\gamma q$ was proved (where $M\rightarrow+\infty$). 
The exponent $\gamma$ was allowed to be any number in the interval $(0,1/2].$ Precisely, Theorem 4.1 of \cite{elsmaily1}, with $\gamma=1/2$, yields that
$$\displaystyle{\lim_{M\rightarrow+\infty}\frac{c_{{\textstyle{\Omega,M\,\textstyle{Id},\sqrt{M}\;q, \,f}}}^*(e)}{\displaystyle{\sqrt{M }}}=2\;\sqrt{f'(0)}.}$$
Using the above together with the continuity of $c^*$ with respect to $M$, there then exists $M_0:=M_0(\delta)>0$ such that
$$\forall\,M\geq M_0(\delta),~~\ds{0<c_{\textstyle \Omega, M\,\textstyle{Id},\sqrt{M}\;q, \,f}^*(e)<\sqrt{M}\left(2\sqrt{f'(0)}+\delta\right)}.$$

 \underline{Step 2.} We fix $M_1>\max(1,M_0(\delta)).$ Then,
\begin{equation}\label{smaily1}
\ds{0<c_{\textstyle{\Omega,\textstyle M_1\,\textstyle{Id},\textstyle \sqrt{M_1}q, \,f}}^*(e)<\sqrt{M_1}\left(2\sqrt{f'(0)}+\delta\right)}.
\end{equation}

\underline{Step 3.} For the fixed number $M_1,$ we also have $\sqrt{M_1}\,q$ in the form of shear flow. We now look at the effect of small diffusion while $\sqrt{M_1}\,q$ is considered as an advection field. Indeed, Theorem 3.3 of El Smaily \cite{elsmaily1} yields that
$$\lim_{\epsilon\rightarrow 0^+}c_{\textstyle{\Omega,\epsilon\,\textstyle{Id},\sqrt{M_1}\;q, f}}^*(e)=\max_{y\in\overline{\omega}}(\sqrt{M_1}q_1(y))=\sqrt{M_1}\;\max_{y\in\overline{\omega}}(q_1(y)).$$
Consequently, there exists $\epsilon_1=\epsilon_1(\delta,M_1)\in (0,M_1)$ (we can choose $\epsilon_1$ as small as we like) such that
\begin{equation}\label{smaily2}
\begin{array}{ll}
\forall \,0<\epsilon\leq \epsilon_1,~   c_{\textstyle{\Omega,\epsilon\,\textstyle{Id},\sqrt{M_1}q, \,f}}^*(e)&>\sqrt{M_1}\ds \max_{y\in\overline{\omega}}(q_1(y))-\delta\\
&>\sqrt{M_1}\left[\ds \max_{y\in\overline{\omega}}(q_1(y))-\delta\right]>0.
\end{array}
\end{equation}

\noindent\underline{Step 4.} Owing to (\ref{smaily}), (\ref{smaily1}) and (\ref{smaily2}),
 one then gets
$$\forall\hbox{ $0<\epsilon\leq\epsilon_1,$}~~c_{\textstyle{\Omega,\epsilon\,\textstyle{Id},\sqrt{M_1}\,q, f}}^*(e)>c_{\textstyle{\Omega,M_1\textstyle{Id},\sqrt{M_1}q, f}}^*(e),$$
and this completes the proof. \hfill$\Box$

\subsection{Case of non-proportional diffusions (answer to the second question)}\label{counter to second question}
In Subsection \ref{counter_to_first_question}, we saw that a prefixed advection field was multiplied by a large enough constant $M_1$ in order to spoil the increasing behavior of the minimal speed with respect to diffusion matrices which are proportional to the identity matrix. A simple scaling of the reaction-advection-diffusion equations that correspond to Theorem \ref{proportional} implies that, under the same notations, 
$$c^*_{\textstyle{\Omega,(\epsilon/ \sqrt{M_1})\textstyle Id,q, f/ \sqrt{M_1} }}(e)>c^*_{\textstyle{\Omega, \sqrt{M_1}\textstyle Id,q,  f/ \sqrt{M_1} }}(e),$$
and by continuity of the minimal speed with respect to diffusion coefficients, one can then find non-proportional matrices $A$ and $B$ (in the neighborhood of $(\epsilon/ \sqrt{M_1})\textstyle Id$ and $\sqrt{M_1}\textstyle Id$ respectively) that satisfy the general assumptions of Section \ref{intro} and such that 
$A<B$ and  $$c^*_{\textstyle{\Omega,A, q, f/\sqrt{M_1} }}(e)>c^*_{\textstyle{\Omega, B,q, f/ \sqrt{M_1} }}(e).$$ This can be summarized as a perturbation argument applied to a strict inequality. However, we notice that the latter inequality corresponds to reaction terms which are small.  The importance of the result we give in this subsection is that it leads to the construction of a counterexample to the monotonicity of the speed with respect to diffusion while $q$ and $f$ are prefixed (neither small nor large).  It turns out that the counterexample given in Theorem \ref{near_0} concerns diffusion matrices which are
 non-proportional by construction and not by a perturbation argument of the result of the previous section. 
\begin{notation}\label{A_b}
For each $b>0,$ let $A_{b}$ denote the $N\times N$ symmetric positive definite matrix
having the form $$A_b=\left(
                        \begin{array}{ccccc}
                          1 & 0 & \ldots & \ldots & 0 \\
                          0& b & \ddots &  & \vdots \\
                          \vdots & \ddots & \ddots & \ddots & \vdots \\
                          \vdots &  & \ddots & \ddots & 0 \\
                          0 & \ldots & \ldots & 0& b \\
                        \end{array}
                      \right)
                         .$$
\end{notation}

\begin{thm}\label{near_0}
Let $e=(1,0,\ldots,0)\in\mathbb{R}^{N},$ $\Omega=\mathbb{R}\times\omega\subseteq\mathbb{R}^{N},$
where $\omega$  can be either bounded or unbounded satisfying (\ref{comega}), and let
$q=(q_{1}(y),0,\ldots,0)$ be a shear flow satisfying (\ref{cq})
where $q_{1}\not\equiv0$ on $\overline{\omega}.$ 
Assume that the nonlinearity satisfies (\ref{kpp_hom_1}-\ref{kpp_hom_2}). For each $b>0,$
consider the reaction-advection-diffusion problem
\begin{eqnarray}\label{A^b with u}
\left\{
  \begin{array}{rl}
    u_t(t,x,y)&=\nabla\cdot(A_{b}\nabla u)+
q_{1}(y)\partial_{x}u+f(u),
\vspace{4 pt}\\
&=\partial_{xx}u+b\Delta_{y}u+q_{1}(y)\partial_{x}u+f(u)\hbox{ in }\mathbb{R}\times\Omega,\vspace{4 pt}\\
 \nu(x,y)\cdot A_{b}\nabla_{x,y}u&=\nu_{\omega}(y)\cdot\nabla_{y}u(t,x,y)=0\hbox{ for }(t,x,y)\in\mathbb{R}\times\mathbb{R}\times\partial\omega,
\end{array}
\right.
\end{eqnarray}
 where $A_{b}$ is the matrix introduced in Notation \ref{A_b}. Then,
\begin{enumerate}
\item \begin{equation}
\lim_{b\rightarrow +\infty}c^*_{\Omega, \textstyle{A_{b}},q,f}(e)=2\sqrt{f '(0)},
\end{equation}
\item \begin{equation}\label{A_b-as_b_goestozero}
 \lim_{b\rightarrow0^+}c^*_{\Omega, \textstyle{A_b},q,f }(e)=\max_{\overline{\omega}}\left(q_1(y)\right)+2\sqrt{f '(0)}.
 \end{equation}
\end{enumerate}
\end{thm}
\textbf{Proof.} Consider the following change of variables
$$\displaystyle{\forall(t,x,y)\in\mathbb{R}\times\mathbb{R}\times\omega,~~
v(t,x,y)=u(t,\frac{x}{\sqrt{b}},y)}.$$ Thus,  $\forall (t,x,y)\in\mathbb{R}\times\mathbb{R}\times\omega,$
$$\displaystyle{v_t(t,x,y)=u_t(t,\frac{x}{\sqrt{b}},y),} \;\displaystyle{\partial_{x}v(t,x,y)=\frac{1}{\sqrt{b}}\partial_x u(t,\frac{x}{\sqrt{b}},y)},$$
$$\displaystyle{\partial_{xx}v(t,x,y)=\frac{1}{b}\partial_{xx}u}(t,\frac{x}{\sqrt{b}},y)\hbox{ and }\Delta_{y}v(t,x,y)=\Delta_{y}u(t,\frac{x}{\sqrt{b}},y).$$
Owing to the invariance of $\Omega$ in the $x-$direction, we have
the boundary condition
$$\forall(t,x,y)\in\mathbb{R}\times\partial\Omega,\quad\nu(x,y)\cdot\nabla_{x,y}v(t,x,y)=0.$$
Consequently, the problem (\ref{A^b with u}) is equivalent to the
problem
\begin{eqnarray}\label{A^M with v}
 \left\{
   \begin{array}{rl}
     v_t(t,x,y)&=\displaystyle{b\partial_{xx}v+b\Delta_{y}v+\sqrt{b}\,q_{1}(y)\partial_{x}v(t,x,y)+f(v),}\vspace{4 pt}\\
     &=\displaystyle{b\,\Delta_{x,y}v+\sqrt{b}\,q_{1}(y)\partial_{x}v+f(v)\hbox{ in }\mathbb{R}\times\mathbb{R}\times\omega,}\vspace{4 pt}\\
\nu(x,y)\cdot\nabla_{x,y}v(t,x,y)&=0\hbox{ for }(t,x,y)\in\mathbb{R}\times\mathbb{R}\times\partial\omega.
   \end{array}
 \right.
\end{eqnarray}

The minimal speed of problem (\ref{A^M with v}) exists in this setting (according to Theorem 1.2 in \cite{bh1}) and is denoted here by $\displaystyle{c^*_{\textstyle{\Omega,b\,\textstyle{Id},{\sqrt{b}}\,q,f}}(e)},$ where $Id$ stands for the identity $N\times N$ matrix. As in the proof of the Theorem \ref{proportional}, we apply Theorem 4.1 of  \cite{elsmaily1} with $\gamma=\displaystyle{{1}/{2}}$ to conclude that
\begin{equation}\label{lim of alpha_M over rad M}
    \displaystyle{\lim_{b\rightarrow+\infty}\frac{c^*_{\textstyle{\Omega,b\,\textstyle{Id},\sqrt{b}q, f}}(e)}{\displaystyle{\sqrt{b }}}=2\sqrt{f'(0)}.}
\end{equation}
On the other hand, looking at the problems (\ref{A^b with u}) and (\ref{A^M with v}) and  the relation between $u$ and $v$ and owing to the minimality (uniqueness) of the threshold $c^*$, one knows that the corresponding minimal speeds have the following relation
\begin{equation}\label{rescale}
\forall b>0,~~\displaystyle{c_{\displaystyle{\Omega,b\,\textstyle Id,\sqrt{b}\;q,
\,f}}^*(e)}=\sqrt{b}\,c_{\Omega, \textstyle{A_{b}},q,f}^*(e).
\end{equation}
Together with (\ref{lim of alpha_M over rad M}), we obtain that
$$\displaystyle{\lim_{b\rightarrow +\infty}c_{\Omega,
\displaystyle{A_{b}},q,f}^*(e)=2\sqrt{f'(0)}}.$$

\textbf{\underline{Proof of (\ref{A_b-as_b_goestozero}).}}
For the limit as $b\rightarrow0^+,$  we know from the rescaling formula (\ref{rescale}) that proving (\ref{A_b-as_b_goestozero}) is equivalent to prove that
\begin{equation}\label{basic}
 \lim_{b\rightarrow0^+}\frac{c_{\textstyle{\Omega,b\textstyle Id,\sqrt{b}\,q, \,f}}^*(e)}{\displaystyle{\sqrt{b }}}=
\max_{\overline{\omega}}\,\left(q_{1}(y)\right)+2\sqrt{f'(0)}.
\end{equation}
Indeed one can apply similar techniques to those used in the proof of Theorem 3.3 of \cite{elsmaily1} ($\zeta$ replaced by  $f'(0)$) to  obtain prove this. However, for the sake of completeness, we do the proof here in details. One knows that $c_{\textstyle{\Omega,b\textstyle Id,\sqrt{b}\,q, \,f}}^*(e)$ is given by the variational formula of \cite{BHN1} 
\begin{equation}\label{varform}
\forall b>0, ~~c_{\textstyle{\Omega,b\textstyle Id,\sqrt{b}\,q, \,f}}^*(e)=\min_{\lambda>0}\frac{k(\lambda)}{\lambda},
\end{equation} 
where $k(\lambda)$ is the principal eigenvalue of the following problem
\begin{equation}\label{eigen_equation_eps_A _with_q}
\begin{array}{lll}
\displaystyle{L_{\Omega,e,b
Id,\,q,f,\lambda}}\,\psi&=&b\Delta\psi+2b\lambda e\cdot\nabla\psi+\sqrt{b}q_{_{1}}(y)\partial_{x}\psi\vspace{4pt}\\
&&+\left[b\lambda^{2}+\lambda \sqrt{b}q_{_{1}}(y)+f'(0)\right]\psi\;\hbox {in }\;\Omega=\mathbb{R}\times\omega,
\end{array}
\end{equation}
with the boundary conditions $\nu \cdot \nabla\psi=\lambda \nu\cdot {e} \psi$ on $\partial\Omega.$ 
The eigenfunction $\psi=\psi(x,y_1,y_2)$ is positive in $\overline{\Omega},$ 1-periodic in $x$, and $(L_1,\ldots,L_{d})-$periodic in $y_1$, and unique up to multiplication by a non-zero constant.
Now, having the coefficient  $q_1$  independent of $x,$ this yields that one can choose $\psi$ independent of $x.$ Hence, the elliptic operator
$\displaystyle{L_{\Omega,e, bId,\,\sqrt{b}q,f,\lambda}}$ can be reduced to the symmetric operator
$$\displaystyle{L_{\Omega,e, b
Id,\,\sqrt{b}q,f,\lambda}}\,\psi=b\Delta\psi+\left[b\lambda^{2}+\lambda\sqrt{b} q_{_{1}}(y)+f'(0)\right]\psi.$$
 Consequently, we have the following Rayleigh quotient over $H^1(C_\omega),$ where $C_\omega$ is the periodicity cell of $\omega$ introduced in Section \ref{intro} above,
\begin{equation}\label{rayleigh}
    \begin{array}{ll}
\forall \lambda>0,\,\forall b>0,\\
\quad-\displaystyle{k(\lambda)}=\displaystyle{\min_{\varphi\in
    H^{1}(C_\omega)\setminus\{0\}}\frac{\displaystyle{b\int_{C_\omega}|\nabla\varphi |^2dy}-\lambda\sqrt{b}\int_{C_\omega}q_{_{1}}(y)\varphi^{2}-\displaystyle{{\int_{C_\omega}\left[{\lambda}^{2}b+f'(0)\right]\varphi^{2}(y)dy}}}{\displaystyle{\int_{C_\omega}\varphi^{2}(y)dy}}.}
\end{array}
\end{equation}
Formula (\ref{rayleigh}) leads directly to the upper bound 
$$\forall\lambda>0,\forall b>0,~~\frac{k(\lambda)}{\lambda\sqrt{b}}\leq \max_{y\in\overline{\omega}}(q_1(y))+\lambda \sqrt{b}+\frac{f'(0)}{\lambda\sqrt{b}}.$$
 Testing the right hand side at $\lambda_0=\sqrt{\frac{f'(0)}{b}},$ one then concludes that 
 $$\forall b>0,~~\min_{\lambda>0}\frac{k(\lambda)}{\lambda\sqrt{b}}\leq \max_{y\in\overline{\omega}}(q_1(y))+2\sqrt{f'(0)},$$
 and hence
 \begin{equation}\label{upper}
\limsup_{b\rightarrow0^+}\frac{c_{\textstyle{\Omega,b\textstyle Id,\sqrt{b}\,q, \,f}}^*(e)}{\sqrt{b}}\leq \max_{y\in\overline{\omega}}(q_1(y))+2\sqrt{f'(0)}.
 \end{equation}
For the sharp lower bound of the $\liminf$ as $b\rightarrow0^+,$ we proceed as follows. Let $y_0$ be a point of $C_\omega\subseteq\overline\omega$ where the periodic function $q_1(y)$ attains its positive maximum (we recall here that $q_1\not\equiv 0,$ $q_1$ is periodic in $y$ and has a zero average over $C_w$ and thus the point $y_0$ is in the interior of the cell $C_\omega$). Then, by the definition of a maximum point of a continuous function, for each $0<\delta< q_1(y_0)$, there exists a neighborhood $U_\delta\subset C_\omega$ of $y_0$ such that 
$$\forall y\in U_\delta,~~q_1(y)\geq q_1(y_0)-\delta\geq0.$$
Now we test the Rayleigh quotient against a smooth function $\phi(y)$ compactly supported in $U_\delta$ and normalized by $\int_{C_\omega}\phi^2dy=1$ to get (after switching from $-k(\lambda)$ to $+k(\lambda)$) 
\begin{equation}\label{lower}
\displaystyle{\frac{k(\lambda)}{\lambda}\geq \lambda b+\displaystyle{\frac{1}{\lambda}}\beta(b)+\sqrt{b}(q_1(y_0)-\delta),}
\end{equation}
where $\beta(b)=f'(0)-b\displaystyle{\displaystyle{\int_{U_\delta}|\nabla\phi|^{2}}}.$ We notice that $\beta(b)>0$ for $0<b<b_0:=\displaystyle{\frac{f'(0)}{\displaystyle{\int_{U_\delta}|\nabla\phi|^{2}}}}$.
The function $\lambda\mapsto  \lambda b+\displaystyle{\frac{1}{\lambda}}\beta(b)$ attains its minimum over $(0,+\infty)$ at $\lambda_*=\sqrt{\frac{\beta(b)}{b}}$ and the value of this minimum is
$2\sqrt{b}\sqrt{\beta(b)}$. Referring to (\ref{lower}), one then gets
$$\forall 0<b<b_0,\forall \lambda>0,~~\frac{k(\lambda)}{\lambda}\geq 2\sqrt{b}\sqrt{\beta(b)}+\sqrt{b}(q_1(y_0)-\delta).$$ The right hand side of the last inequality is independent of $\lambda.$ Thus, the  variational formula (\ref{varform}) yields that
$$\forall 0<b<b_0,~~\frac{c_{\textstyle{\Omega,b\textstyle Id,\sqrt{b}\,q, \,f}}^*(e)}{\sqrt{b}}\geq 2\sqrt{\beta(b)}+(q_1(y_0)-\delta),$$ for all $0<\delta<q_1(y_0).$
 Thanks to the non-dependence of $U_\delta$ on $b,$ we pass to the limit as $b\rightarrow0^+$ to get
 $$\forall 0<\delta<q_1(y_0), ~~\liminf_{b\rightarrow0^+}\frac{c_{\textstyle{\Omega,b\textstyle Id,\sqrt{b}\,q, \,f}}^*(e)}{\sqrt{b}}\geq 2\sqrt{f'(0)}+(q_1(y_0)-\delta).
 $$   
 This is enough to conclude that
 \begin{equation}\label{liminf}
 \liminf_{b\rightarrow0^+}\frac{c_{\textstyle{\Omega,b\textstyle Id,\sqrt{b}\,q, \,f}}^*(e)}{\sqrt{b}}\geq 2\sqrt{f'(0)}+(q_1(y_0))
 \end{equation}
and, together with (\ref{upper}), finish the proof of (\ref{A_b-as_b_goestozero}).
\hfill$\Box$

 \subsubsection*{A counterexample as an application of Theorem \ref{near_0}}
 \vskip 0.2cm
\noindent Let $e=(1,0,\ldots,0)$ and
$\Omega=\mathbb{R}\times\omega$ (in particular $\Omega$ can be taken as the whole space $\R^N$).  Choose $f=f(u),$ and $q=(q_{_{1}}(y),0,\cdots,0)$ with
$\displaystyle{\int_{C}q_{_{1}}(y)dy=0}$ and $q_1\not\equiv 0.$ Thus there exists
$\delta>0$ such that 
$$\displaystyle{2\,\sqrt{f^{'}(0)}+\delta<\max_{y\in\overline{\omega}}\left(q_{_{1}}(y)\right)+2\,\sqrt{f^{'}(0)}-\delta.}$$
It follows, from Theorem \ref{near_0} above, and from the continuity of the speed with respect to the diffusion factor (variational formula of \cite{BHN1})
that there exist $\varepsilon_{0}>0$ and $M_{0}>0$ such that
$$\forall 0<\varepsilon\leq\varepsilon_{0},\quad c_{\Omega,\textstyle{A_{\varepsilon}},q,f}^*(e)>\max_{y\in\overline{\omega}}\left(q_{_{1}}(y)\right)+2\,\sqrt{f '(0)}-\delta\quad\hbox{and}$$
$$\forall \,M\geq M_{0}>0,\quad c_{\Omega, \textstyle{A_{M}},q,f}^*(e)<2\sqrt{f^{'}(0)}+\delta.$$

Consequently, choosing $\varepsilon$ small enough and $M$ large
enough, it follows that $\displaystyle{A_M\geq A_{\varepsilon}}$ in
the sense of order relation on positive definite matrices; however,
$$c_{\Omega, \textstyle{A_{M}},q,f}^*(e)<c_{\Omega, \textstyle{A_{\varepsilon}},q,f}^*(e).$$

Therefore the answer to the second question is negative in general 
even when the advection is a fixed shear flow and the nonlinearity $f$ is fixed and homogenous.\hfill $\Box$

\section*{Acknowledgments}
I would like to thank Professor Fran\c{c}ois Hamel for bringing this problem to my attention during my Ph.D work under his supervision. Also, I would like to thank him for his valuable, fruitful and precise suggestions while reading a preprint of this work.

\end{document}